\documentclass[preprint,12pt]{elsarticle}

\usepackage{amsthm,amsmath,amssymb}

\usepackage[]{algorithm2e}

\journal{Journal of Discrete Algorithms}

\begin{document}

\begin{frontmatter}

\title{Fast Algorithm for Enumerating Diagonal Latin Squares of Small Order}
\author[ISDCT]{Stepan Kochemazov\corref{cor1}}
\ead{veinamond@gmail.com}
\author[KSU]{Eduard Vatutin}
\ead{evatutin@rambler.ru}
\author[ISDCT]{Oleg Zaikin}
\ead{zaikin.icc@gmail.com}
\cortext[cor1]{Corresponding author}

\address[ISDCT]{Matrosov Institute for System Dynamics and Control Theory SB RAS\\Lermontova str. 134, Irkutsk, Russia, 664033}
\address[KSU]{Southwest State University\\ 50 let Octyabrya str. 94, Kursk, Russia, 305040\\}

 


\begin{abstract}
In this paper we propose an algorithm for enumerating diagonal Latin squares of small order. It relies on specific properties of diagonal Latin squares to employ symmetry breaking techniques, and on several heuristic optimizations and bit arithmetic techniques to make use of computational power of state-of-the-art CPUs. Using this approach we enumerated diagonal Latin squares of order at most 9, and vertically symmetric diagonal Latin squares of order at most 10.   
\end{abstract}
\end{frontmatter}

\section{Introduction}\label{sec:intro}
A Latin square of order $N$ is an $N\times N$ table, filled with elements from the set $\{0,\ldots,N-1\}$ so that in each row and column each element appears exactly once \cite{Colbourn2006}. If in addition its main diagonal and main antidiagonal both contain every possible element from $0$ to $N-1$ then a Latin square is called diagonal. We call a Latin square vertically symmetric if a sum of two elements in the same row, positioned symmetrically relative to the virtual vertical middle line is equal to $N-1$.
What makes Latin squares interesting is that they represent a good example of a very well studied combinatorial design with numerous applications, for which nevertheless there exist several exceptionally hard related open problems. Probably the most well known one is to determine if there exists the set of three mutually orthogonal Latin squares of order $10$ \cite{DBLP:journals/moc/EganW16}.
Due to combinatorial nature of Latin squares, there naturally arise various enumeration problems, classification problems, etc. Usually, it is considered unrealistic to perform explicit enumeration of all Latin squares of a specific order, however, the development of computers and algorithms made it possible to obtain several results of such kind in recent years, in particular, to enumerate Latin squares of orders 10 and 11 \cite{McKay1995,McKay2005}. However, as far as we are aware, the number of diagonal Latin squares even for order 8 was not known until the end of 2016. That is why in the present paper we develop an approach aimed at enumeration of diagonal Latin squares of small order.

In the context of enumeration, the main difference between diagonal Latin squares and ordinary Latin squares  consists in the fact that diagonal Latin squares form much smaller equivalence classes, because the uniqueness constraints on diagonals essentially restrict most transformations used to form equivalent Latin squares (mainly, arbitrary row and column permutations). 
In this paper we first design a fast algorithm for explicit generation of diagonal Latin squares, augment it with the so-called symmetry breaking techniques that make partial use of equivalence classes of diagonal Latin squares, and apply it to enumerate diagonal and vertically symmetric diagonal Latin squares of small order. 

There are two main contributions in the present paper. The first one is the optimized brute-force algorithm for enumeration of diagonal Latin squares and related designs, such as Latin rectangles, etc. Essentially, the algorithm represents a Latin square as an integer array and uses $\leq N^2$ nested loops to traverse all possible variants of Latin square cells values. Its simple structure can be improved by several heuristic-based optimizations. In particular, the order, in which the cells are filled, greatly influences the algorithm performance. Also, the implementation details, i.e. how we perform necessary checks and assignments, how we organize each loop, etc., play an important role as well. The bit arithmetic techniques greatly aid the performance of these operations. The resulting version of the algorithm makes it possible to enumerate up to 7 millions of diagonal Latin squares of order 9 per second on one CPU core. The symmetry breaking techniques that are based on the class of transformations which convert diagonal Latin squares into diagonal Latin squares, allow to reduce the size of the search space by several orders of magnitude. We used the constructed algorithm to obtain the second contribution, i.e. to enumerate diagonal Latin squares of order at most 9 and vertically symmetric diagonal Latin squares of order 10 (the corresponding numbers were unknown before). 

Let us present a brief outline of the paper. In the next section we discuss possible ways to generate Latin squares. Then we describe the basic structure of our algorithm that we use as a basis of further optimizations. In Section \ref{sec:bit} we show how bit arithmetic techniques make it possible to greatly increase the performance of the algorithm in practice, and experimentally evaluate different algorithm versions. In Section \ref{sec:eq} we describe how the equivalence classes of diagonal Latin squares can be constructed and use this information to introduce into the proposed algorithm symmetry breaking techniques. Then we describe our computational experiments, in the course of which we enumerated diagonal Latin squares of order at most 9, enumerated vertically symmetric diagonal Latin squares of order 10, and estimated the number of diagonal Latin squares of order 10. After this we discuss related works and draw conclusions.

\section{Algorithm Description}\label{sec:alg}
Hereinafter, without the loss of generality we assume that generation and enumeration mean the same and treat these terms as interchangeable. Within the context of enumeration it is sensible to consider only algorithms that are deterministic and complete, i.e. the ones that can generate all possible representatives of the desired species which satisfy fixed constraints. Since we do not intend to store generated diagonal Latin squares, the enumeration should proceed in a fixed order and should not employ randomization on any stage, so that we process the whole search space, and at the same time do not enumerate some diagonal Latin square more than once.

In the next subsection let us consider several algorithmic concepts that fit the description above. Since our main goal is to enumerate diagonal Latin squares of order 9, we mainly consider and evaluate possible algorithms in application to this problem. If not stated otherwise, all performance evaluations are performed on one core of Intel Core i7-6770 CPU, 16 Gb RAM. All algorithms proposed in the paper were implemented in C++ using Microsoft Visual Studio 2015 compiler for Windows or \texttt{gcc} (different versions) for Linux. In pseudocode we will mostly follow the C++ syntax.

\subsection{Approaches to Generating Diagonal Latin Squares}\label{subsec:ls_gen}

Each row and each column of a Latin square is a permutation of $N$ elements. It means that for small $N$ one can generate all possible permutations and construct Latin squares by combining them. For example, it is possible to fill the square row by row, meanwhile checking that different rows do not have equal elements in the same positions. However, in this case once several rows are filled, the number of available variants for the remaining rows drops very significantly, thus making simple exhaustive search, that tries to put every possible permutation as the next row, very ineffective.
For example, if we consider Latin squares of order 10, we have $10!= 3\;628\;800$ possible permutations. We can put each of them as the first row, then for the second row we loop through the list and test if a permutation number $i$ does not violate Latin square constraints. For rows after the 5th the number of such permutations (that can be put as the next row) is in the range of hundreds at most. Thus if we cycle through all available permutations to put into, say, 8th row, even if we can test if they fit or not very fast - the process is quite ineffective.

In this context it is sensible to represent the original problem as \textit{exact cover} instance  \cite{DBLP:conf/coco/Karp72} and employ relatively sophisticated algorithms, such as DLX \cite{DLX}, which can restrict the search space "on-the-fly". Note, that if we are interested only in diagonal Latin squares -- the case presents two more uniqueness constraints, that need to be taken into account. Our preliminary evaluations showed that bit arithmetic-aided exhaustive search and DLX make it possible to enumerate about $5\times 10^5$ diagonal Latin squares of order 9 per second on one CPU core.

In the present paper we follow another simple approach to generating Latin squares. Within it we represent Latin square of order $N$ as an array of $N^2$ integer values, corresponding to its cells, and fill their values in a fixed order. In the most basic variant we implement this enumeration procedure in the form of $N^2$ nested \texttt{for} loops. 
On the first glance, it seems that this approach is too crude and should lose compared to the ones mentioned above. Indeed, if we fill square elements from left to right from top to bottom, then the generation speed is very low: about $6\times 10^3$ diagonal Latin squares of order 9 per second on one CPU core. However, after several optimizations, that we consider below, this approach significantly outperforms the others.

\subsection{Algorithm Design}\label{subsec:alg_design}
Assume that we consider enumeration of diagonal Latin squares of order $N$.
For this purpose our algorithm uses several auxiliary constructs:
\begin{enumerate}
	\item Integer array $LS[N][N]$ which contains a Latin square.
	\item Integer arrays $Rows[N][N]$ and $Columns[N][N]$, where we reflect which elements are already "occupied" in each row/column.
	\item Integer arrays $MD[N]$ and $AD[N]$ where we reflect which elements are "occupied" on main diagonal and main antidiagonal.
	\item Integer value $SquaresCnt$ in which we accumulate the number of squares.
\end{enumerate}

{
\begin{algorithm}[htbp]
\label{alg1}
 \KwData{$LS[N][N]$, $Rows[N][N]$, $Columns[N][N]$, $MD[N]$, $AD[N]$, $SquaresCnt$}
 \tcc{All variables are initialized by 0.} 
\tcc{Iterate over all possible values of cell [0][0]}
\For {$LS[0][0]=0$; $LS[0][0]<N$; $LS[0][0]=LS[0][0]+1$}{
	\tcc{If value is occupied within the row, column or diagonals - continue to the next value of loop variable}
	\lIf{$Rows[0][LS[0][0]] || Columns[0][LS[0][0]] || MD[LS[0][0]]$}	{continue}
	\tcc{Otherwise mark the value as occupied and proceed}
	$Rows[0][LS[0][0]]=1$\;
	$Columns[0][LS[0][0]]=1$\;
	$MD[LS[0][0]]=1$\;
	
	\For {$LS[0][1]=0$; $LS[0][1]<N$; $LS[0][1]=LS[0][1]+1$}{
			\lIf{$Rows[0][LS[0][1]] || Columns[1][LS[0][1]]$}	{continue}			
				$Rows[0][LS[0][1]]=1$\;
				$Columns[1][LS[0][1]]=1$\;
				
				...\\
				
				\tcc{Increment SquaresCnt if reached the last element}
				\For {$LS[N-1][N-1]=0$; $LS[N-1][N-1]<N$; $LS[N-1][N-1]=LS[N-1][N-1]+1$}{
				\lIf{$Rows[N-1][LS[N-1][N-1]] || Columns[N-1][LS[N-1][N-1]] || MD[LS[N-1][N-1]]$}	{continue}
				$SquaresCnt=SquaresCnt+1$\;
				}				
				...\\												
				
				$Rows[0][LS[0][1]]=0$\;
				$Columns[1][LS[0][1]]=0$\;				 	
		}
	
	\tcc{On exit from the loop mark value as `free'}	
	$Rows[0][LS[0][0]]=0$\;
	$Columns[0][LS[0][0]]=0$\;
	$MD[LS[0][0]]=0$\;
}
\caption{General outline of the algorithm}
\end{algorithm}
}

As we mentioned above, the order, in which we fill cells greatly influences the performance of the algorithm. However, for now, let us introduce the general outline of the algorithm for enumerating diagonal Latin squares with simple order when we fill the square from the first (topmost leftmost) element to the last. Its pseudocode is presented as Algorithm \ref{alg1}.

Note, that we can already make one simple optimization. It is clear that each non-diagonal Latin square can be effectively transformed (by means of row and column permutations) to a Latin square, in which the first row and the first column appear in ascending order $0,1,\ldots,N-1$ (the corresponding procedure is usually referred to as normalization \cite{Colbourn2006}). It means that we can safely fix the values of corresponding variables in the array $LS[N][N]$ and modify the initialization stage. As a result we have $(N-1)^2$ inner loops instead of $N^2$. For diagonal Latin squares we fix either the first row, the first column or the main diagonal to $0,1,\ldots,N-1$, thus having $N^2-N$ inner loops. Hereinafter we assume that the first row is fixed in all algorithms and experiments.

To make further constructions easier it is natural to represent the algorithm as a sequence of loops. The order in which the cells values are filled is then reflected by the order of these loops. The structure of a loop is presented as Algorithm \ref{alg2}.

\begin{algorithm}[htbp]
\label{alg2}
 \KwData{$LS[N][N]$, $Rows[N][N]$, $Columns[N][N]$, $MD[N]$, $AD[N]$, $SquaresCnt$, $i$, $j$}
\tcc{Iterate over all possible values of cell [i][j]}
\For {$LS[i][j]=0$; $LS[i][j]<N$; $LS[i][j]=LS[i][j]+1$}{
	\tcc{Check if the value is occupied in the current row, column or diagonals}    
    \tcc{Without the loss of generality, we assume that the element lies on the intersection of Main diagonal and Main antidiagonal. If it does not, then we omit corresponding entries ($MD[LS[i][j]]$ and/or $AD[LS[i][j]]$ when computing the value of $Condition_{i,j}$)}
    bool $Condition_{i,j} = Rows[i][LS[i][j]] || Columns[j][LS[i][j]] || MD[LS[i][j]] || AD[LS[i][j]]$\;
    \uIf{$Condition_{i,j}$}{
    \tcc{If $Condition_{i,j}=True$ it means that the value LS[i][j] is already occupied. In this case we proceed to next value of LS[i][j]}    
    continue\;}
	\tcc{Otherwise mark the value as occupied and proceed}
	$Rows[i][LS[i][j]]=1$\;
	$Columns[j][LS[i][j]]=1$\;
	\tcc{Important: entries for diagonals are included only if $i=j$ and/or $i+j=N-1$}
	$MD[LS[i][j]]=1$\;
    $AD[LS[i][j]]=1$\;	
	\textbf{BODY OF INNER LOOP FOR NEXT CELL $LS[i'][j']$\;}	
	$Rows[i][LS[i][j]]=0$\;
	$Columns[j][LS[i][j]]=0$\;
	\tcc{Important: entries for diagonals are included only if $i=j$ and/or $i+j=N-1$}
    $MD[LS[i][j]]=0$\;
    $AD[LS[i][j]]=0$\;	
}
\caption{Inner loop structure}
\end{algorithm}

So within inner loop we cycle over all possible values from $0$ to $N-1$ to put into cell $LS[i][j]$. We use auxiliary arrays to store information whether the value $l$ was already used within row/column/main diagonal/main antidiagonal. In particular, $Rows[i][l]=1$ if and only if the value $l$ is assigned to some cell within the $i$-th row, and $Columns[j][l]=1$ if the value $l$ is assigned to some cell within the $j$-th column. The same goes for $MD[l]=1$ and $AD[l]=1$ for main diagonal and main antidiagonal. Once we find the value that can be put into $LS[i][j]$ without violating any constraint, we do it and refresh the information about occupied values within $Rows[i]$, $Columns[j]$ and arrays for diagonals $MD$, $AD$, if applicable. After this we proceed to the inner loop for next cell $LS[i'][j']$. On exit from the loop, we clear values corresponding to $LS[i][j]$ in $Rows[i]$, $Columns[j]$, $MD$ and $AD$ to prepare for the next iteration.

In the body of the loop for the last cell according to specified order we only increment the counter $SquaresCnt$ (since if we reached it, it means that we successfully constructed a Latin square).  Now let us consider the question of the optimal order of cells.

\subsection{On the Optimal Order of Cells}\label{subsec:cells_order}

In our experiments we noticed a very interesting pattern. It turned out, that the algorithm performance, especially, if we consider diagonal Latin squares, greatly depends on the order, in which the cells are filled. In particular, when we change only the order of cells, and do not touch any other parameters, the average generation speed may vary from several thousand to several hundred thousand diagonal Latin squares of order 9 per second. After detailed empirical evaluation, we figured the strategy that works best. In essence, it implements the ideas suggested in \cite{Golomb:1965:BP:321296.321300}, that in the backtrack search one should narrow the search space as much as possible on each step.

Let us now consider our implementation of the strategy that yields the order, with which the proposed algorithm shows the best performance. Let us consider how it works for diagonal Latin squares of order 9. In accordance with the general outline of the algorithm, we fill the cells of a Latin square $LS=\{LS[i][j]\}$. The first row of this square is fixed: $LS[0][j]=j$, $j=0,\ldots,8$. We use the iterative process to choose the cell to be assigned next. It is important to note, that in this process we use only the information whether the cell is already assigned, and do not know exact value stored by any cell. It is easy to see, that each cell $LS[i][j]$ is involved in at least two and at most four "uniqueness constraints": one for the $i$-th row, one for the $j$-th column, and two more for the main diagonal and main antidiagonal if $i = j$ and/or $i=8-j$. Let us consider the value $V_{i,j}^k=r_i^k+c_j^k+md_k(i,j)+ad_k(i,j)$. Here $k$ is the iteration number, $r_i^k$ is the number of assigned cells in the $i$-th row on the $k$-th step, $c_j^k$ -- the number of assigned cells in the $j$-th column on the $k$-th step, $md_k(i,j)$ is the number of assigned cells on the main diagonal if $i=j$ and $0$ otherwise, and $ad_k(i,j)$ is the number of assigned cells on the main antidiagonal if $i=8-j$, and $0$, otherwise. In a sense, the number $V_{i,j}^k$ reflects how "constrained" is a cell with indexes $i,j$. On each iteration step we choose the cell that has the largest $V_{i,j}^k$. If several cells have the same value of $V$, we choose the first of them when ordered in lexicographic order.

It also makes sense to use an additional simple heuristics: if after choosing a new cell at some step $k'$ the number of assigned cells in some row, column, main diagonal or main antidiagonal becomes $N-1$, then the remaining cell is automatically assigned next because we can compute it directly thanks to corresponding uniqueness constraints. Now let us return to the case of diagonal Latin squares of order 9. On the $0$-th iteration only the cells in the first row are assigned. It means that the most constrained (in the aforementioned sense) cell is the one that lies on the intersection of main diagonal and main antidiagonal, it has $V_{4,4}=3$, while for all other cells it is at most $2$. Then we proceed as outlined above. As a result, we obtain the order of cells presented in Figure  \ref{fig:dls9_cells_order}.

\begin{figure}[htbp]
\begin{center}
\begin{tabular}{ccccccccc}
-&-&-&-&-&-&-&-&-\\
20&2&16&17&21&18&19&3&22\\
25&26&6&23&27&24&7&28&29\\
55&56&57&10&53&11&58&59&60\\
61&63&65&45&1&47&67&69&70\\
62&64&66&12&54&14&68&71&72\\
32&33&8&30&34&31&9&35&36\\
39&4&42&37&40&38&43&5&41\\
13&49&50&44&48&46&51&52&15
\end{tabular}
\end{center}
\caption{Order of cells for generation of diagonal Latin squares of order 9 (the first row is fixed a priori, so it is omitted)}
\label{fig:dls9_cells_order}
\end{figure}

In essence, as a result of the outlined procedure for diagonal Latin squares of order 9 we start with diagonal elements, and then fill the rest. When we embed this order into the algorithm, it can enumerate about 1.2 million diagonal Latin squares of order 9 per second on one CPU core. It is interesting, that when applied to constructing the optimal order for ordinary Latin squares, the proposed heuristics constructs the trivial order of cells: row by row, column by column. Let us now proceed to other optimizations.

\subsection{Optimizations}\label{subsec:optim}
The following two techniques, while quite simple, make it possible to increase the performance of the above algorithm to  about 1.8 million squares per second.

\subsubsection{Use formula to compute the last element in a row/column/diagonal}\label{subsubsec:formula}
At certain point within the algorithm, there appear situations, when in some row, column or main diagonal/antidiagonal, there are $N-1$ assigned cells out of $N$. Due to the fact that these elements are the subject of "uniqueness" constraints, in these cases it is possible to compute the value of a remaining cell directly, thus eliminating the need to introduce a loop for this purpose. Without the loss of generality, let us assume that we assigned the first $N-1$ of $N$ elements in the $j$-th row. Then the formula for the remaining element looks as follows:
\begin{equation*}
LS[j,N-1]=N\times(N-1)/2-\sum\limits_{l=0}^{N-2}LS[l,j].
\end{equation*}
Of course, we need to make sure that the obtained value does not violate other uniqueness constraints before proceeding deeper into the search space.

\subsubsection{Lookahead heuristic}\label{subsubsec:lookahead}
We borrowed the next technique from the area of combinatorial search. It represents a kind of a lookahead heuristic \cite{Dechter2003117}. The basic idea is that on some levels of the search (i.e. in some of the loops) there arise situations, when as a result of assigning value to a current cell the amount of constraints on some other cells within the same row/column/diagonal may exceed $N$, thus there is a possibility that there are no possible assignments for these "over-constrained" cells. So we can spend some resources and look ahead before branching further, spending a little more computational resources now, so that we avoid spending much more of them later.

For simplicity, assume that for non-diagonal element $LS[i][j]$ the following situation arises (remind that $Rows[i][l]=1$ if and only if the value $l$ is assigned to some cell within the $i$-th row, and $Columns[j][l]=1$ if the value $l$ is assigned to some cell within the $j$-th column):
\begin{equation*}
\sum\limits_{l=0}^{N-1} (Rows[i][l] \vee Columns[j][l]) = N.
\end{equation*}
It is clear, that in this case we can stop looking further, since the currently examined portion of the search space is reduced to empty set. Thus we revert the last assignment of Latin square cell and proceed.

It is important, that this heuristic should be used with care. Conditional operators in large quantities can easily slow the search down, making any performance gain disappear. Therefore, there should be found a tradeoff, by choosing on which levels to apply it. After empirical evaluation and testing, for enumeration of diagonal Latin squares of order 9 we determined, that it works best when we apply lookahead within inner loops from number 51 to 60.

Now let us consider how to apply bit arithmetic techniques to improve the algorithm performance further.

\section{Bit Arithmetic Implementation}\label{sec:bit}
In order to improve the performance of the suggested algorithm we can do several things: merge/remove repeated actions, do the same things faster and  reduce the number of conditional operators involved. In this section we show that it is possible to do all this thanks to employing bit arithmetic techniques. Hereinafter, without the loss of generality, we assume that all integer values contain at least $16$ bits. 

Let us describe the modifications we introduce to the algorithm outlined above in order to transition to using bit arithmetic. First, we drop one dimension from arrays $Rows$, $Columns$, $MD$, $AD$ since we fuse it within one integer value with $\geq N$ bits. Second, we represent the values of Latin square cells in a different manner: instead of $LS[i][j]=k$ we now write $LS[i][j]=1<<k$, where $<<k$ means left bit shift for $k$ positions. We also introduce an array of auxiliary variables $CR[N][N]$ to keep track of the number of current constraints on each Latin square element.

The modified inner loop structure is presented as Algorithm \ref{alg3}. without the loss of generality we assume that the considered element lies on the intersection of main diagonal and main antidiagonal. In case it does not, the corresponding entries are simply removed (i.e. for the most simple case $CR[i][j]=Rows[i]|Columns[j]$) together with operators marking the new value occupied/free in $MD$ and $AD$.

\begin{algorithm}[ht]
\label{alg3}
 \KwData{$LS[N][N]$, $CR[N][N]$, $Rows[N]$, $Columns[N]$, $MD$, $AD$, $SquaresCnt$,$i$,$j$}
\tcc{Compute vector of possible values for cell [i][j]}
$CR[i][j]= Rows[i] | Columns[j]|MD|AD$\;
	\tcc{Iterate over all possible values of cell [i][j]}
\For {$LS[i][j]=1$; $LS[i][j]<(1<<N)$; $LS[i][j]=LS[i][j]<<1$}{
	\tcc{Check if the value is occupied }
	\lIf{$(CR[i][j]\&LS[i][j])!=0$}	{	continue	}
	\tcc{Otherwise mark the value as occupied and proceed}
	$Rows[i]=Rows[i]|LS[i][j]$\;
	$Columns[j]=Columns[j] | LS[i][j]$\;
	$MD=MD| LS[i][j]$\;
    $AD=AD| LS[i][j]$\;
	\textbf{BODY OF INNER LOOP FOR NEXT CELL $LS[i'][j']$}\;	
	$Rows[i]=Rows[i]\oplus LS[i][j]$\;
	$Columns[j]=Columns[j] \oplus LS[i][j]$\;
    $MD=MD \oplus LS[i][j]$\;
    $AD=AD \oplus LS[i][j]$\;
	}
\caption{Inner loop structure with bit arithmetic}
\end{algorithm}

Here the main performance gain is thanks to the use of an array $CR[i][j]$. By definition it contains $1$-bits in all positions in which $LS[i][j]$ \textbf{can not} take the value of 1. Therefore, we can obtain at once the spectrum of available values, instead of checking the availability of each cell value by iterating over them. When implemented in the proposed manner, the algorithm is able to generate about $2.6\times 10^6$ diagonal Latin squares of order 9 per second even without employing optimizations from the end of the previous section. Nevertheless, there is still a room for improvement. Despite the fact that we now compute the vector of possible values for each particular element only once, we still need to iterate over all $N$ values of $LS[i][j]$. Is there a way to make proper use of this information? 

Fortunately, yes. In particular, we can reconstruct the \textbf{for} loop in order to iterate over only such values of $LS[i][j]$ that satisfy the condition $CR[i][j]\&LS[i][j]=0$, thus eliminating the need for \textbf{if} block in the inner loop body in Algorithm \ref{alg3}. For this purpose we introduce a new constant $AllN=1<<(N-1)$ that has exactly $N$ $1$-bits in the beginning.  The next version of the algorithm heavily relies on the bit twiddling tricks that make it possible to isolate the rightmost 1-bit ($y=x\&(-x)$) and to turn off the rightmost 1-bit ($y=x \& (x-1)$). Let us present the modified inner loop structure in the pseudocode as Algorithm \ref{alg4}.

\begin{algorithm}[ht]
\label{alg4}
 \KwData{$LS[N][N]$, $CR[N][N]$, $L[N][N]$, $Rows[N]$, $Columns[N]$, $MD$, $AD$, $SquaresCnt$,$i$,$j$}
\tcc{Compute vector of possible values for cell [i][j]}
$CR[i][j]= Rows[i] | Columns[j] | MD | AD$\;
	\tcc{Iterate over values of cell [i][j] that do not violate any uniqueness constraint.}
\For {$L[i][j]=AllN\oplus CR[i][j]$; $L[i][j]!=0$; $L[i][j]=L[i][j]\&(L[i][j]-1)$}{
	$LS[i][j]=L[i][j]\&(-L[i][j])$\;
	\tcc{Mark the value as occupied and proceed}
	$Rows[i]=Rows[i]|LS[i][j]$\;
	$Columns[j]=Columns[j] | LS[i][j]$\;
	$MD=MD| LS[i][j]$\;
    $AD=AD| LS[i][j]$\;
	\textbf{BODY OF INNER LOOP FOR NEXT CELL $LS[i'][j']$}\;	
	$Rows[i]=Rows[i]\oplus LS[i][j]$\;
	$Columns[j]=Columns[j] \oplus LS[i][j]$\;
    $MD=MD  \oplus LS[i][j]$\;
    $AD=AD \oplus LS[i][j]$\;
	}
\caption{Optimized inner loop structure with the bit arithmetic}
\end{algorithm}

The first major achievement in the improved inner loop design is that we got rid of one of the two conditional operators (the remaining one is incorporated into the \texttt{for} cycle construction and it is unlikely that we can get rid of it). So, we first compute the value of $CR[i][j]$. It contains the bit vector with $1$ bits in positions corresponding to values of $LS[i][j]$ that violate any of uniqueness constraints. Then we use additional auxiliary integer array $L[N][N]$. In the \texttt{for} cycle we initialize $L[i][j]$ with possible values of $LS[i][j]$ that \textbf{do not} violate any constraint and iterate over them by switching off the rightmost $1$ bit until $L[i][j]$ becomes $0$. For each value of $L[i][j]$ we produce the value of $LS[i][j]$ by isolating the rightmost $1$-bit in $L[i][j]$. Once $L[i][j]$ becomes $0$ it means that we processed all available alternatives. This improved algorithm version makes it possible to generate about $6\times 10^6$ diagonal Latin squares of order 9 per second without heuristic optimizations.

Note,  that we can use the BMI1 instruction set, supported in the state-of-the-art CPUs, via the corresponding 
intrinsics\footnote{https://software.intel.com/sites/landingpage/IntrinsicsGuide/} to switch off the rightmost $1$ bit (\texttt{\_blsr\_u64}) and isolate the rightmost $1$ bit (\texttt{\_blsi\_u64}) faster. On average, it leads to a performance increase of 5-10\%.

\begin{table}[htbp]
	\caption{Performance of the proposed versions of the algorithm for generation of Latin squares of small order.}
	\centering
		\begin{tabular}{|c|c|c|}
		\hline
		\textbf{Version}&\textbf{Problem}&\textbf{Squares per second}\\
		\hline
		
		Standard (Algorithms \ref{alg1} and \ref{alg2}) &DLS9&$1.8\times 10^6$\\
		\hline
		
		Bit arithmetic (Algorithm \ref{alg3})&DLS9&$2.6\times 10^6$\\
		\hline
		
		&DLS9&$6.8\times 10^6$\\
		\cline{2-3}
		&LS8& $9\times 10^6$\\
		\cline{2-3}
		Optimized bit arithmetic  &DLS8 & $5.8 \times 10^6$\\
		\cline{2-3}
		(Algorithm \ref{alg4})&LS9 & $8.0 \times 10^6$\\
		\cline{2-3}
		&LS10 & $6.3\times 10^6$\\
		\cline{2-3}
		&DLS10 & $ 6.0\times 10^6$\\
		\hline
		\end{tabular}	
	\label{tab:versions_performance}
\end{table}

In Table \ref{tab:versions_performance} we show the generation speed for different classes of Latin squares. Here for diagonal Latin squares we  fixed the first row (in an ascending order) and for ordinary Latin squares we fixed both the first row and the first column (in an ascending order) because if we know the number of normalized (diagonal) Latin squares -- it is easy to compute the corresponding number of non-normalized squares. Table entry $(D)LS$ followed by number stands for (diagonal) Latin squares of specific order. The order of cells in each case was determined according to heuristic procedure outlined in Subsection \ref{subsec:cells_order}. The performance of Algorithm 4 for $DLS9$ was measured for the algorithm versions that use Lookahead heuristic (in other cases the corresponding optimization requires a lot of empirical evaluation and testing, so for other entries the lookahead heuristic is not used). 

It is clear that bit arithmetic techniques make it possible to significantly increase the performance of the algorithm. Note, that nor here, nor in the next section we did not use processor intrinsics to boost the performance because we mainly employed the computers and clusters that did not support BMI1 instruction set.

\section {Equivalence Classes for Diagonal Latin squares} \label{sec:eq}

The so-called symmetry breaking techniques \cite{Walsh2006} are often used in algorithms dealing with combinatorial objects in order to avoid visiting the parts of the search space that are in some way similar to that already processed. In the case of Latin squares when all possible species representatives are to be processed, their space is usually divided into the so-called main classes and then only one representative for each class is considered. This question is covered in more detail, for example, in \cite{Colbourn2006,McKay2005}. Informally speaking, a main class is mainly formed by Latin squares that can be produced from each other by all possible transpositions of rows and columns, since these transformations do not lead to violation of any constraints. However, the diagonal Latin squares represent a special case, because in their case most transpositions lead to violation of uniqueness constraints on  elements from the main diagonal and the main antidiagonal. Below we present the class of symmetric row-column transpositions which transform diagonal Latin square to a diagonal Latin square.

\subsection{Transformations That Preserve Elements on Diagonals}
 The roots of the corresponding class of transformations lie in the area of Magic squares \cite{10.2307/j.ctt7s1bz}, thus the transformations themselves are called \textit{M-transformations}. A magic square is an $n\times n$ table filled with integer numbers in such a way that the sums of numbers in each row, each column and also in the main diagonal and main antidiagonal are all equal to the same so-called ``magic constant''. It is clear that each diagonal Latin square is a Magic square. To the best of our knowledge, the M-transformations, while widely known, have not been published in any one paper to be cited, thus we will briefly present them below. 
There are three types of basic M-transformations.
\begin{itemize}
\item The first type is formed by mirroring a diagonal Latin square horizontally (or vertically) and relative to the main diagonal (or antidiagonal). There are 4 variants of this transformations.

\item The second type contains all transpositions of at least two columns that are positioned symmetrically with respect to the middle with simultaneous transposition of two symmetrically positioned rows, for example, transposition of $0$-th and $(n-1)$-th columns with simultaneous transposition of $0$-th and $(n-1)$-th row. The number of such transformations for diagonal Latin squares of order $n$ is $2^{\lfloor \frac{n}{2} \rfloor}$ (equal to the number of all subsets of a set with $\lfloor\frac{n}{2}\rfloor$ elements).

\item In the third type there are all transpositions where we simultaneously transpose at least two columns in the left half of a Latin square  and at least two columns positioned symmetrically with respect to the middle in the right half of a square with simultaneous similar transposition of rows, for example transposition of $0$-th column with $1$-th column, $(n-2)$-th column with $(n-1)$-th column, $0$-th row with $1$-th row and $(n-2)$-th row with $(n-1)$-th row. The number of these transformations for diagonal Latin square of order $n$ is $\lfloor \frac{n}{2}\rfloor !$ (equal to the number of permutations of a set with $\lfloor \frac{n}{2}\rfloor$ elements).
\end{itemize}

The intuition behind basic M-transformations is relatively simple: each time we transpose rows or columns, we do it in such a way, that the elements from both main diagonal and main antidiagonal, positioned within these rows or columns, are only transposed within their diagonal. Because of this reason, any combination of basic M-transformations, to which we refer as M-transformation, also preserves elements on diagonals. Hereinafter we assume that after each M-transformation we fix the first row of a constructed diagonal Latin square to ascending order by renaming elements.

As a result, we can form an equivalence class for (normalized) diagonal Latin squares of order $N$ of size at most $4 \times 2^{\lfloor \frac{N}{2}\rfloor} \times \lfloor \frac{N}{2}\rfloor!$. Meanwhile for ordinary Latin squares a main class can contain up to $6\times N\times N!$ (normalized) Latin squares. The difference is drastic: for $N=9$ it is 1 536 normalized diagonal Latin squares vs 19 595 520. For $N=10$ it is 15 360 vs 217 728 000. 

It is currently not clear if it is possible to generate/enumerate only unique representatives of equivalence classes of diagonal Latin squares without significantly sacrificing the generation performance, especially if we take into account the fact that in the context of enumeration for each equivalence class representative we have to find the power of its equivalence class. We found that for our purposes of enumeration of diagonal Latin squares of small order there is a relatively easy workaround that makes it possible to mostly preserve the achieved algorithm effectiveness. 

\subsection{Symmetry Breaking}
The interesting observation consists in the following: we can apply M-transformations to partial diagonal Latin squares. Assume that we have an incomplete diagonal Latin square, in which only elements from the first and last row, and from the main diagonal and main antidiagonal are present. Let us informally refer to such an incomplete diagonal Latin square as to \textit{hourglass} design. In Figure \ref{hourglass} we show an example of hourglass design for $N=10$.
\begin{figure}[htbp]
\begin{center}
\begin{tabular}{cc}
$
\left(
\begin{array}{cccccccccc}
 0& 1& 2& 3& 4& 5& 6& 7& 8& 9\\
\_& 4&\_&\_&\_&\_&\_&\_& 6&\_\\
\_&\_& 1&\_&\_&\_&\_& 3&\_&\_\\
\_&\_&\_& 6&\_&\_& 7&\_&\_&\_\\
\_&\_&\_&\_& 7& 1&\_&\_&\_&\_\\
\_&\_&\_&\_& 2& 8&\_&\_&\_&\_\\
\_&\_&\_& 8&\_&\_& 2&\_&\_&\_\\
\_&\_& 0&\_&\_&\_&\_& 9&\_&\_\\
\_& 5&\_&\_&\_&\_&\_&\_& 3&\_\\
 4& 2& 3& 1& 8& 6& 9& 0& 7& 5\\
\end{array}
\right)
$
&
$
\left(
\begin{array}{cccccccccc}
 9& 8& 7& 6& 5& 4& 3& 2& 1& 0\\
\_& 6&\_&\_&\_&\_&\_&\_& 4&\_\\
\_&\_& 3&\_&\_&\_&\_& 1&\_&\_\\
\_&\_&\_& 7&\_&\_& 6&\_&\_&\_\\
\_&\_&\_&\_& 1& 7&\_&\_&\_&\_\\
\_&\_&\_&\_& 8& 2&\_&\_&\_&\_\\
\_&\_&\_& 2&\_&\_& 8&\_&\_&\_\\
\_&\_& 9&\_&\_&\_&\_& 0&\_&\_\\
\_& 3&\_&\_&\_&\_&\_&\_& 5&\_\\
 5& 7& 0& 9& 6& 8& 1& 3& 2& 4\\
\end{array}
\right)
$
\end{tabular}
\end{center}
\caption{Example of hourglass design for $N=10$. The design on the right is a vertically mirrored version of the one on the left.}
\label{hourglass}
\end{figure}
How do we apply M-transformations to an incomplete diagonal Latin square? The same way we do it to an ordinary diagonal Latin square, however, we only set values to elements which are known. Similarly, when renaming elements to fix the first row in ascending order, we rename only known elements. For example, if we mirror the hourglass design from the left part of Figure \ref{hourglass} vertically, we will produce a hourglass design from the right part of Figure \ref{hourglass}.

It is easy to see, that M-transformations, that do not transpose the first row with any other row besides the last one, and do not mirror the square relative to its diagonals, transform  hourglass design into hourglass design. Assume that $H_1$ and $H_2$ are two distinct hourglass designs of order $N$, and $H_2$ is produced from $H_1$ by applying to it some M-transformation, to which we refer to as $\mu$, i.e. $H_2=\mu(H_1)$.

\textbf{Proposition 1.} The number of diagonal Latin squares of order $N$ that share $H_1$ is equal to the number of diagonal Latin squares of order $N$ that share $H_2=\mu(H_1)$. 

\textbf{Sketch Proof.}

When we apply $\mu$ to a hourglass design $H_1$, we implicitly apply it to every single diagonal Latin square which shares cell values with $H_1$. Since $\mu(H_1)=H_2$ it means that any diagonal Latin square that shares cell values with $H_1$ will be transformed to a diagonal Latin square that shares cell values with $H_2$.

An interesting and quite obvious property of M-transformations is that for each transformation there is an inverse transformation. It means that for any combination of basic M-transformations $\mu$ there is a combination of basic M-transformations (which cancel them) $\mu^{-1}$ such that for an arbitrary diagonal Latin square $A$ the following holds $\mu(\mu^{-1}(A))=\mu^{-1}(\mu(A))=A$. That is why, the situation, when there are two distinct Latin squares $A$ and $B$ that share $H_1$ and $\mu(A)=\mu(B)$, is impossible.

From the above it follows that the sets of diagonal Latin squares that share $H_1$ and $H_2$ have the same power. $\qed$

\textbf{Corollary 1.}
It is possible to transform the proposed enumeration algorithm for diagonal Latin squares of order $N$ in the following way:
\begin{enumerate}
\item Construct all possible hourglass designs of order $N$.
\item Split the space of hourglass designs of order $N$ into equivalence classes using applicable M-transformation. Compute the power of each equivalence class.
\item Process each constructed equivalence class as follows.
\begin{enumerate}
\item Choose one representative of an equivalence class and enumerate all diagonal Latin squares of order $N$ that share cell values with this representative.
\item Multiply the result of enumeration by the power of the equivalence class.
\end{enumerate}
\item Summarize the results.
\end{enumerate}
\subsection{Embedding Symmetry Breaking into the Algorithm}
The proposed symmetry breaking technique fits well into the general idea of our algorithm: we fill Latin square cells in specific order, so we can with little effort fill in hourglass design elements before everything else. For this purpose we do not even need to change the constructed optimal cells order that much: the first row is fixed beforehand, and the  diagonals are filled before everything else. Thus to produce hourglass designs we need only to fill the last row of cells right after filling diagonals. 

\begin{algorithm}[htbp]
\SetKwFunction{FCanonize}{Canonize}
\SetKwProg{Fn}{Function}{:}{}
\Fn{\FCanonize{$Hourglass$}}{
\KwIn{$Hourglass[N][N]$, $TR_0[2]$, $TR_1[2^{\lfloor \frac{N}{2}\rfloor}$, $TR2[(\lfloor\frac{N}{2}\rfloor-1)!]$} 
\tcc{Hourglass[N][N] -- Hourglass design}
\tcc{$TR_0[2]$, $TR_1[2^{\lfloor \frac{N}{2}\rfloor-1}]$, $TR_2[(\frac{N}{2}-1)!]$ -- arrays of M-transformations of three types.}
\tcc{EQC - container for constructed hourglass designs}
int $EQC[2 \times 2^{\lfloor \frac{N}{2}\rfloor} \times (\lfloor \frac{N}{2}\rfloor-1)!][N][N]$)\;
\tcc{K - counter for constructed hourglass designs}
int $K=0$\;
\tcc{F -- flag indicating if $Hourglass$ is canonic form}
bool $F=True$\;
\tcc{Iterate over equivalent hourglass designs.}
\For {$i_1=0;\;i_1<2;\;i_1=i_1+1$}{
  \For {$i_2=0;\;i_2<2^{\lfloor\frac{N}{2}\rfloor};\;i_2=i_2+1$}{
    \For {$i_3=0;\;i_3<(\lfloor\frac{N}{2}\rfloor-1)!;\;i_3=i_3+1$}{		
        int $T[N][N]$\tcp*{Apply M-transformations.}		
        $T=Transform(Hourglass,TR_0[i_1], TR_1[i_2], TR_2[i_3])$\;
        \tcc{Normalize constructed hourglass design.}		
        $Normalize(T)$\;
        \uIf{$T<Hourglass$}{
        F=false\;
        break\;
        }
        \uElse{
         $K=K+1$\;
         $EQC[K]=T$\;
        }
    }
  }
}
\uIf{$F=True$}{
\tcc{Compute the number of distinct hourglass designs in EQC.}		
int $R=ComputeDistinct(EQClass,K)$\;
\Return{R}\;
}
\uElse{
\Return{0}\;
}
}
\label{alg5}
\caption{Algorithm for computing canonic form of hourglass design}
\end{algorithm}

However, things stop being simple at this point because the number of hourglass designs can be very large. Indeed, for diagonal Latin squares of order $N=8$ it is equal to 22 192 248. One equivalence class for hourglass designs of order $8$ contains at most 192 designs, so as a result we still need to store about 100 000 (116 857) classes representatives. For larger dimensions the number of hourglass designs is far too large to store it in a memory of a state-of-the-art computer. 

However, there is a computational way to choose unique class representatives, which, while presenting overhead in terms of computational resources, uses very little memory. We can embed equivalence checking right into the general algorithm. Hereinafter we assume that a \textit{canonic form} of a hourglass design is a hourglass design from its equivalence class which is lexicographically the first for a fixed ordering of elements (for example, in natural order from left to right from top to bottom). It means that once we generated a hourglass design $H$, we can apply to it the procedure, that determines if it is equal to its canonic form. This procedure applies to $H$ all applicable  M-transformations and if some equivalent design $H'$ is lexicographically less, $H'<H$, then it means that $H$ is not a canonic form. In this case the procedure terminates and returns zero. Otherwise, if the procedure applied all possible combinations of M-transformations and $H$ is still the least of them, then we know exactly the power of $H$'s equivalence class, and the procedure returns this number. Note, that in practice for different reasons, the hourglass designs $H_1$, $H_2$ produced by applying different M-transformations $\mu_1$, $\mu_2$, $\mu_1 \neq \mu_2$ to hourglass design $H$ can be equal ($H_1=\mu_1(H)=\mu_2(H)=H_2$), thus we in fact need to store all constructed hourglass designs and, in case our $H$ is a canonic form, to compute the number of distinct designs in the equivalence class. The pseudocode of this algorithm is presented as Algorithm \ref {alg5}.

Now that we introduced the function to canonize hourglass designs, we can embed it into the general algorithm. Informally speaking, we split the algorithm into two parts: the first (outer loops) that generates a hourglass design, and the second (inner loops) that enumerates all possible diagonal Latin squares that share generated hourglass design if it coincides with its canonic form. The general outline of the corresponding algorithm is presented as Algorithm \ref{alg6}

\begin{algorithm}[ht]
\KwData{$LS[N][N]$, $Rows[N]$, $Columns[N]$, $MD$, $AD$, $SquaresCnt$,$TotalCnt$} 
$TotalCnt=0$\;
\For {...}{
	\tcc{Inner loops generating hourglass design}
	
	\tcc{Once a hourglass design is generated check if it is equal to canonic form}    
    int $multiple=Canonize(LS)$\;
    SquaresCnt=0\;
    \uIf {$multiple\neq 0$}{
    	\For {...}{
        	\tcc{If a hourglass design is equal to canonic form, then we enumerate all diagonal Latin squares that share it, otherwise we generate next hourglass design}
        }    
    }   	
    $TotalCnt = TotalCnt + SquaresCnt \times multiple$\;
        	...\\
\tcc{ends of inner loops generating hourglass designs}
}	
\label{alg6}
\caption{General outline of enumeration algorithm with symmetry breaking}
\end{algorithm}

In total the proposed scheme makes it possible to reduce the search space for $N=8$ and $N=9$ by almost 200 times.

\section{Computational Experiments}
We applied the proposed algorithm to solve three problems: enumeration of diagonal Latin squares of orders 8 and 9, enumeration of vertically symmetric diagonal Latin squares of order 10, and estimation of the number of diagonal Latin squares of order 10. Note, that we first applied the algorithm in its pure form (without symmetry breaking techniques) to enumerate diagonal Latin squares of order 9 in two large scale computational experiments and to estimate the number of diagonal Latin squares of order 10. Later, using new version of the algorithm we verified the achieved result and enumerated vertically symmetric diagonal Latin squares of order 10.

\subsection{Enumeration of diagonal Latin squares of order 8}\label{sec:exp8}
We used the preliminary version of the presented algorithm (in the form outlined in Section \ref{sec:alg}, i.e. without bit arithmetic-related optimizations) to enumerate diagonal Latin squares of order 8 (7 447 587 840) \cite{VatutinZZMKT2016}, which is reflected in the sequence A274171 \cite{A274171} in the On-Line Encyclopedia of Integer Sequences. At the time of experiment it took about 30 CPU hours. Our current algorithm implementation achieves this result in 21 minutes on one CPU core. If we augment it with symmetry breaking techniques, this time is even less - about 90 seconds on one CPU core. Theoretical speedup must be about 200 times, but in practice it is only 14. The explanation lies in the fact that symmetry breaking presents computational overhead, and for simple problems its impact is significant. If we exclude the time required to split hourglass designs of order $8$ into equivalence classes, then the actual enumeration process takes about 6 seconds.

\subsection{Enumeration of diagonal Latin squares of order 9}\label{sec:exp}
At first we performed two separate experiments to enumerate diagonal Latin squares of order 9. The first one was held in the volunteer computing project Gerasim@home \cite{VatutinVT2015}, while the second was performed within the computing cluster ``Academician V.M. Matrosov'' of Irkutsk supercomputing center SB RAS\footnote{http://www.hpc.icc.ru}. 

In both cases we decomposed the problem in the following way. Note, that we fix the first Latin square row in an ascending order. Since we fill in the cells of Latin square in a specific order (presented at Fig. \ref{fig:dls9_cells_order}), it means that we can choose small number of the first cells to be filled and process their correct assignments separately. In our experiments we used in the role of "decomposition set" the first 10 cells (according to the aforementioned order). There are 1 225 884 possible correct assignments of these cells. We constructed a workunit for each of these assignments. In volunteer computing a workunit is a batch of computational tasks. In each workunit the problem was to enumerate all possible diagonal Latin squares of order 9 with the fixed values of the cells in the first row and 10 more cells (in total 19 cells fixed out of 81). 

Since the proposed workunits do not depend on each other in any way, it means that we can process them according to the embarrassing parallelism concept \cite{Foster:1995:DBP:527029}. In particular, it allows us to employ both parallel and distributed computing systems for this purpose. As a result we obtain an array of 1 225 884 integer numbers. Their total sum is equal to the number of diagonal Latin squares of order 9 with fixed first row. 

\subsubsection{Experiment in a volunteer computing project}\label{subsec:volunteer}
\begin{figure}[htbp]
	\centering
		\includegraphics[width=11cm]{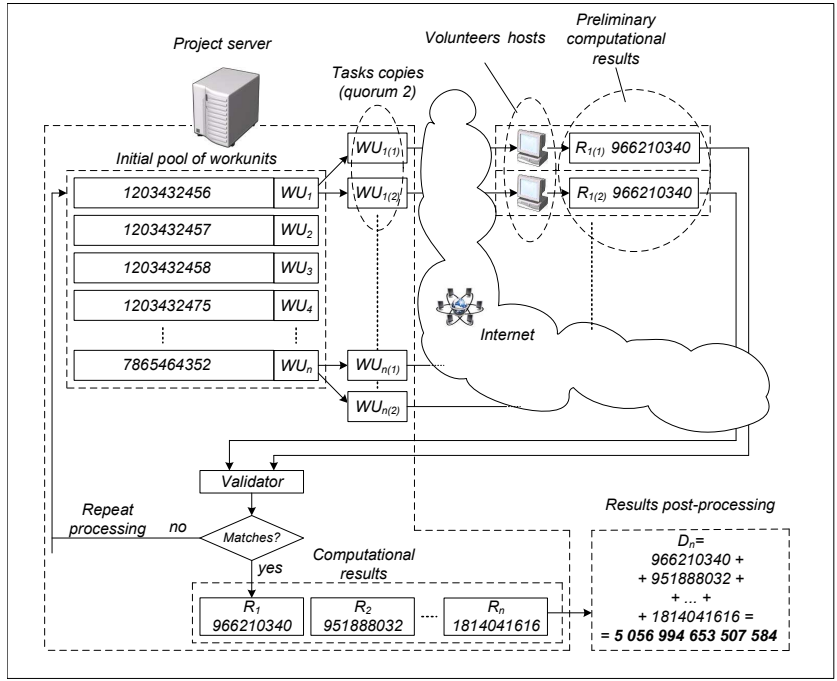}
	\caption{Scheme of the Gerasim@home experiment}
	\label{fig:Gerasim_experiment_scheme}
\end{figure}

The volunteer computing project Gerasim@home \cite{VatutinVT2015} is based on BOINC (Berkeley Open Infrastructure for Network Computing \cite{DBLP:conf/ccgrid/AndersonF06}) platform. One of the goals of this experiment was to prove that volunteer computing can be used to compute the number of specific combinatorial designs because usually, there are relatively high demands to such experiments.
The scheme of this experiment is showed at Fig. \ref{fig:Gerasim_experiment_scheme}.

Let us comment on this scheme. There is a project server, which distributes workunits and processes their results. On each host (usually a personal computer of a volunteer) the client application processes the input data. In our case as an input data the application takes the file containing the workunit description in the form of values of specified 10 cells (as it was outlined above). The number of diagonal Latin squares, obtained as a result of enumeration process, is written to an output file and sent back to server. In order to decrease the influence of possible software and hardware errors we used the so called quorum of 2. It means that, according to the BOINC redundancy technique, two copies of each workunit are sent to different volunteers. After processing the results are compared and if they coincide, then the result of the corresponding workunit is marked as correct. If not, then two new copies of the workunit are generated.

The computational experiment aimed at enumeration of diagonal Latin squares of order 9 was launched on 18 June, 2016. At that moment the performance of the project was about 2.5 teraflops. In the experiment we used only x86 and x64 client applications for Windows OS. The experiment took 3 months and ended on 17 September, 2016. In total, about 500 volunteers from 51 countries participated in this experiment. They connected to the project about 1000 hosts with the peak performance of 5 teraflops (an average performance was about 3 teraflops).

As a result of post-processing, we determined the number of diagonal Latin squares of order 9 with first row fixed in ascending order: 5 056 994 653 507 584. If we multiply it to $9!$ we obtain the number of diagonal Latin squares of order 9. It should be noted, that if we didn't improve the algorithm performance by means, described in Sections \ref{sec:alg} and \ref{sec:bit}, the corresponding experiment in Gerasim@home would take about 10 years.

\subsubsection{Experiment on a Computing Cluster}\label{subsec:cluster}
In application to hard enumeration problems it is crucial to cross-check the results, since it is possible that small errors remain undetected in specific circumstances and ruin the correctness. One might say that it is especially so when using volunteer computing, however, our empirical results prove otherwise.
In any case, we decided to check ourselves and launched one more experiment aimed at solving the same problem. This experiment was performed within the computing cluster ``Academician V.M. Matrosov'' of Irkutsk supercomputing center SB RAS. Each node of this cluster is equipped with two 16-core AMD Opteron 6276 CPUs and 64 gigabytes of RAM. We used the same approach to decomposition as in Gerasim@home experiment. However, we did not use any redundant calculations, so each workunit was processed exactly once.

We developed an MPI-program (here MPI stands for Message Passing Interface) based on our algorithm.  In this program one process is a control process, and all the remaining processes are computing processes. The control process creates and maintains the pool of workunits to be processed by computing processes. It also accumulates and processes their results. Overall, the organization of the experiment was quite similar to that in Gerasim@home (see Fig. \ref{fig:Gerasim_experiment_scheme}).

The experiment was launched on 17 July, 2016. It took several launches of the MPI-program to finish it. In these launches the number of employed cluster nodes varied from 10 to 15, and the duration varied from 2 hours to 7 days. The majority of launches used 15 nodes with a duration of 7 days. The experiment ended on 17 October, 2016. On average, the experiment took 2 months of computing time with 15 nodes. As a result we verified that the number computed in Gerasim@home was correct. 

\subsubsection{Experiment With Symmetry Breaking}\label{subsec:exp_sym9}
We developed the symmetry breaking technique after two experiments described above already ended, but decided to check the correctness of the result and the technique. We performed the experiment with symmetry breaking on a new part of the computing cluster with performance of about 15 previously used cluster nodes. The experiment ended after six and a half hours and proved the result obtained before.

\subsection{Estimation of the Number of Diagonal Latin Squares of Order 10}\label{subsec:exp_DLS10}
After enumerating diagonal Latin squares of order 9, we naturally decided to apply our algorithm to diagonal Latin squares of order 10. However, it became clear quite fast, that their number is very large. To estimate it we employed the Monte Carlo method \cite{Metropolis49} in the following form. If we specify some order, in which we fill the cells of a diagonal Latin square, we, in fact, can consider an incomplete diagonal Latin square, formed by the first $k$ cells according to the specified order. It is natural to consider the trivial order: the first $k$ elements of a Latin square from left to right from top to bottom. Let us refer to such incomplete diagonal Latin squares of order 10 as to $DLS_{10}^k$. First, for a specific $k$ we compute the number of possible $DLS_{10}^k$, to which we refer as $N_{10}^k$. Then we form a random sample of $DLS_{10}^k$. For each incomplete diagonal Latin square from the sample we  enumerate all possible diagonal Latin squares of order 10 that can be constructed by filling unassigned cells of this $DLS_{10}^k$. As a result of processing of the random sample we construct an estimation of an expected value of the number of diagonal Latin squares of order 10 that share the same $DLS_{10}^k$. By multiplying this estimation of expected value to $N_{10}^k$ we construct the estimation of the number of diagonal Latin squares of order 10.

First we need to choose $k$ in such a way that we can compute the number of $DLS_{10}^k$ in reasonable time, and then be able to process the sample of $DLS_{10}^k$ of a sufficient size. We fix the elements of the first row of a Latin square in an ascending order for simplicity. We started from value $k=30$. The corresponding $N_{10}^{30}$ is 284 086 571 712. However, for each $DLS_{10}^{30}$ it takes several days on one core of state-of-the-art CPU to enumerate all possible diagonal Latin squares that share some $DLS_{10}^{30}$. Thus we chose $k=32$ and computed the estimation for this value. The number of corresponding incomplete diagonal Latin squares $N_{10}^{32}$ is 12 611 543 636 160. We generated a random sample of size 10 000 $DLS_{10}^{32}$ instances and used it to estimate the expected value of the number of diagonal Latin squares of order 10 with fixed $DLS_{10}^{32}$. The corresponding expected value was equal to 11 931 268 344. Thus the estimated number of diagonal Latin squares of order 10 is about $1.5\times 10^{23}$. 

\subsection{Enumeration of Vertically Symmetric Diagonal Latin Squares of Order 10} \label{subsec:exp_sym_SDLS10}
Since the estimated number of diagonal Latin squares of order $10$ turned out to be too large to enumerate them in affordable time given available resources, even with symmetry breaking techniques, we decided to apply our approach to enumeration of a much smaller class of diagonal Latin squares of order 10, in particular, of vertically symmetric diagonal Latin squares. 
A vertically symmetric diagonal Latin square is a diagonal Latin square for which the following holds: 
$$\forall i,j\in\{0,\ldots,N-1\} LS[i][j] = N-1-LS[i][N-1-j]$$
It is clear, that the algorithm for enumerating such squares is much more simple than that for DLS of order $10$: we basically can fill only the left half of the Latin square, because the values of elements in the right half are derived explicitly from the left one. Our algorithm when modified to take vertical symmetry into account makes it possible to enumerate up to 30 million diagonal Latin squares of order $10$ per second on one CPU core.

Note, that M-transformations described above actually preserve vertical symmetry. It means that we can apply the proposed symmetry breaking techniques to this problem as well. The equivalence class for vertically symmetric hourglass designs of order $10$ contains up to $2^5 \times 4! = 768$ entries because mirroring vertically is neutralized by vertical symmetry. Thanks to this it was possible to enumerate vertically symmetric diagonal Latin squares of order 10 quite fast: it took 735 seconds when launched on 8 threads of Intel Core i7-6770 CPU. The resulting number is 82 731 715 264 512\footnote{https://oeis.org/A287649}. 

\section{Related Work}\label{sec:related}
Authors are not aware of algorithms developed specifically for enumeration of diagonal Latin squares. Papers \cite{Bammel:1975:N99:2625518.2625683,McKay1995,McKay2005} describe the approaches that led to enumeration of Latin squares of orders 9, 10 and 11. The corresponding algorithms heavily rely on the ability to permute rows and columns to construct equivalence classes and evaluate their properties. To the best of our knowledge the results of these papers are not applicable to diagonal Latin squares because the vast majority of row-column permutations break diagonal property. Also, 
 diagonal Latin squares form relatively small equivalence classes.

There are several examples of application of parallel and volunteer computing to the search for combinatorial designs based on Latin squares. With the help of a computing cluster there was proven that there is no finite projective plane of order 10 \cite{LAM}. In the volunteer computing project SAT@home several dozen pairs of mutually orthogonal diagonal Latin squares were found \cite{DBLP:conf/mipro/ZaikinKS16}. 

Quite similar approach to the one used in our paper was employed in \cite{DBLP:journals/em/McGuireTC14}. In that paper the hypothesis about the minimal number of clues in Sudoku was proven. The authors developed the fast algorithm to enumerate and check all possible Sudoku variants. The algorithm was implemented and launched on a modern computing cluster. It took about 11 months for this cluster to check all variants. The volunteer computing project Sudoku@vtaiwan \cite{Lin:2010:SMS:1934910.1935592} was used to confirm the solution of this problem.

\section{Conclusions and Future Work}\label{sec:concl}
In this paper we presented the fast algorithm for enumeration of diagonal Latin squares of small order. We developed the fast implementation of the algorithm that exploits the features of state-of-the-art CPUs. Using the proposed symmetry breaking techniques it is possible to significantly reduce the search space. In future we plan to study if the proposed algorithms can be extended to GPUs and used to solve other open problems in related areas.

The present article is a significantly reworked and extended variant of the paper \cite{Vatutin2017-PCT}. The modifications include (but not limited to) new section on symmetry breaking and applications of the corresponding techniques (Section \ref{sec:eq}), thanks to which the main result of \cite{Vatutin2017-PCT} was improved by approximately 200 times, and the related experiments described in subsections \ref{sec:exp8}, \ref{subsec:exp_sym9},  \ref{subsec:exp_DLS10}, \ref{subsec:exp_sym_SDLS10}). 

\section*{Acknowledgements}\label{sec:ack}
We are grateful to Ian Wanless for clarifications regarding the applicability of the algorithm from \cite{McKay2005} to enumeration of diagonal Latin squares, to Alexey Belishev for information about M-transformations, to Alexey Zhuravlev for proposing to use the algorithm with nested loops for enumeration of diagonal Latin squares. We thank all Gerasim@home volunteers and enthusiasts, whose computers took part in the experiment.

Stepan Kochemazov and Oleg Zaikin are supported by Russian Science Foundation (project No. 16-11-10046) and additionally supported by Council for Grants of the President of the Russian Federation (stipends SP-1829.2016.5 and SP-1184.2015.5, respectively).
Eduard Vatutin is partially supported by Council for Grants of the President of the Russian Federation (grant MK-9445.2016.8) and by Russian Foundation for Basic Research (grant 17-07-00317-a).

\textbf{Author contribution.} S. Kochemazov: bit arithmetic implementation of algorithm outlined in Section \ref{sec:bit}, Symmetry breaking techniques and their implementation, described in Section \ref{sec:eq}; 
E. Vatutin: general scheme of algorithm in Section \ref{sec:alg}, the experiment described in Subsection \ref{subsec:volunteer};
O. Zaikin: The experiments described in Subsections \ref{subsec:cluster}, \ref{subsec:exp_sym9}, \ref{subsec:exp_sym_SDLS10}, \ref{subsec:exp_DLS10}.

\bibliographystyle{elsarticle-num} 
\bibliography{jda}

\end{document}